%% file: paper.tex
\documentclass[12pt]{article}

\usepackage{color,amsmath,amssymb,epsf,epsfig,chicago,psfrag,afterpage}
\usepackage{subfigure,subfig,rotating,amsthm}
\usepackage{a4,latexsym}
\usepackage{setspace,calc,longtable,lscape}

\include{mymaths}

\include{theorems}

\newcommand{\figh}{2.5}
\newcommand{\figw}{2.5}

\begin{document}

\title{A New Class of Skewed Bimodal Distributions}

\author{Ricardo S. Ehlers\thanks{Corresponding author. Email: ehlers@icmc.usp.br} \\ 
  {\it Department of Applied Mathematics and Statistics}\\ 
  {\it University of S\~ao Paulo - Brazil}}

\date{}

\maketitle

\section{Introduction}\label{sec:intro}

Probability distributions that can acommodate the possible presence of
heavy tails and skewness in the distribution of a phenomenon have been
the focus of interest in recent years. See for example, \citeN{azz85},
\citeN{fsteel98}, \citeN{azz03}, \citeN{jones03} and \citeN{fsteel06} to
name but a few. However, these distributions fail to capture a
possible bimodality in the data under study. In this paper, our aim is
to introduce a new family of distributions that is flexible 
enough to support skewness, heavy tail and bimodal shapes.

Recently, \citeN{olivero-etal} and \citeN{rocha-etal} extended
Azzalini's skew normal family of distributions to accommodate such
behaviour in the resulting distribution. These authors propose to
disturb the symmetry of the density,
\begin{equation}\label{eq:g}
g(x) = \left(\frac{1+\alpha x^2}{1+\alpha b}\right) f(x), ~~x\in\mathbb{R},
\end{equation}
where $\alpha\ge 0$ and $f(\cdot)$ is symmetric and unimodal with finite
second moment $b$. The parameter $\alpha$ control the uni or
bimodality of $g(\cdot)$ since the density is unimodal if 
$\alpha\in [0,0.5)$ and bimodal if $\alpha\ge 0.5$. Then, they use a
cumulative distribution function $H(\cdot)$ as a skewing mechanism and the
proposed skewed version (possibly bimodal) is given by,
$$
s(x|\alpha,\lambda,H)=
2\left(\frac{1+\alpha x^2}{1+\alpha b}\right) f(x)H(\lambda x),
$$
where the parameter $\lambda\in\mathbb{R}$ introduces skewness.

In this paper, we propose a different route. We first obtain the
skewed version of a unimodal symmetric density using a skewing
mechanism that is not based on a cumulative distribution
function. Then we disturb the unimodality of the resulting skewed
density using the same mechanism as in (\ref{eq:g}). In order to
introduce skewness we use the general method
proposed in \citeN{fsteel98} which transforms
any continuous unimodal and symmetric distribution into a skewed one
by changing the scale at each side of the mode. 
They proposed the following
class of skewed distributions indexed by a shape parameter $\gamma>0$,
which describes the degree of asymmetry,
\begin{equation}\label{eq:skew}
s(x|\gamma)=\frac{2}{\gamma+1/\gamma}
\left\{
  f\left(\frac{x}{\gamma}\right)I_{[0,\infty)}(x)+ 
  f(x\gamma)I_{(-\infty,0)}(x)
\right\},
\end{equation}
where $f(\cdot)$ is a univariate density symmetric around zero and
$I_C(\cdot)$ is an indicator function on $C$.
Note that $\gamma=1$ yields the symmetric distribution as
$s(x|\gamma=1)=f(x)$.
Right skewness corresponds to $\gamma>1$ while left skewness
corresponds to $\gamma<1$.
Our preference for this skewing mechanism is mainly due to its
simplicity and generality. 
Moments calculation is
straightforward if the moments of the underlying symmetric
distribution are available and it does not require calculation of
cumulative distribution functions, which yields faster computations.
Also, it entirely separates the effects of
the skewness and tail parameters thus making prior independence
between the two a plausible
assumption, and hence facilitates the choice of their prior
distributions.

{\proposition \label{p1}
  Let $f$ be a symmetric unimodal density with mode
  zero. If $s(\cdot|\gamma)$ is as defined in (\ref{eq:skew})
  and $b_{\gamma}=\int_{-\infty}^{\infty}x^2s(x|\gamma)dx<\infty$ then,
\begin{eqnarray*}\label{eq:s1}
s(x|\alpha,\gamma) &=&
\left(\frac{1+\alpha x^2}{1+\alpha b_{\gamma}}\right)
\frac{2}{\gamma+1/\gamma}
\left\{
  f\left(\frac{x}{\gamma}\right)I_{[0,\infty)}(x)+ 
  f(x\gamma)I_{(-\infty,0)}(x)
\right\}\\
&=&
\left(\frac{1+\alpha x^2}{1+\alpha b_{\gamma}}\right)
s(x|\gamma)
\end{eqnarray*}
is a density for any $\alpha\ge 0$ and $\gamma>0$.
}

\noindent {\bf Proof.} Clearly $s(x|\alpha,\gamma)\ge 0$. Also,
\begin{eqnarray*}
\int_{-\infty}^{\infty}
\left(\frac{1+\alpha x^2}{1+\alpha b_{\gamma}}\right)s(x|\gamma)dx
&=&
\frac{1}{1+\alpha b_{\gamma}}
\left[
  1+\alpha\int_{-\infty}^{\infty}x^2 s(x|\gamma)dx
\right] = 1
\end{eqnarray*}
since the integral on the right hand side is simply
$b_{\gamma}$. \qed 

\vskip 0.5cm

The existence of
the moments of (\ref{eq:skew}) depends only on the existence
of moments of the symmetric density $f(\cdot)$
and does not depend on $\gamma$. The $r$th moment is given by,
$$
E(X^r|\gamma) = 
\frac{\gamma^{r+1} + (-1)^{r}/\gamma^{r+1}}{\gamma+1/\gamma} ~m_r,
$$
where
$$
m_r = 2\int_0^{\infty}x^r f(x)dx
$$
is the $r$-th absolute moment of $f(x)$ on the positive real
line. It is not difficult to see that when the original symmetric
distribution has mean zero and variance one then $m_2=1$. In this case,
the second moment $b_{\gamma}$ is given by,
$$
b_{\gamma}=\frac{\gamma^3+1/\gamma^3}{\gamma+1/\gamma}.
$$

\noindent So, the moments of this bimodal skewed distribution are given by,
$$
E(X^r|\alpha,\gamma) = 
\frac{1}{1+\alpha b_{\gamma}}
~[E(X^r|\gamma)+\alpha E(X^{r+2}|\gamma)].
$$

\begin{figure}[p]\centering
\begin{tabular}{cc}
\subfigure[][$\gamma=1.5$]
   {\epsfig{file=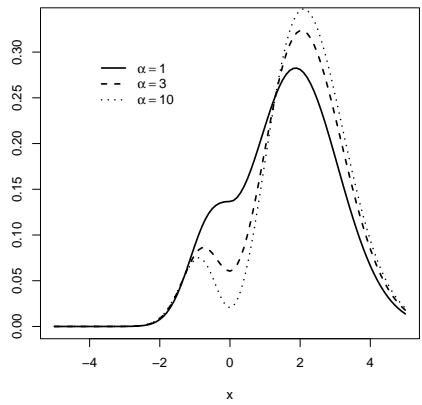,height=\figh in,width=\figw in}}& 
\subfigure[][$\gamma=0.5$]
   {\epsfig{file=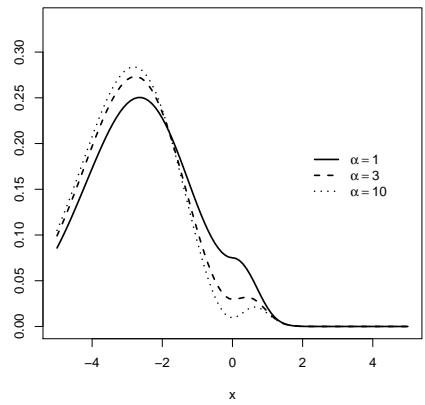,height=\figh in,width=\figw in}}\\
\subfigure[][$\gamma=1.1$]
   {\epsfig{file=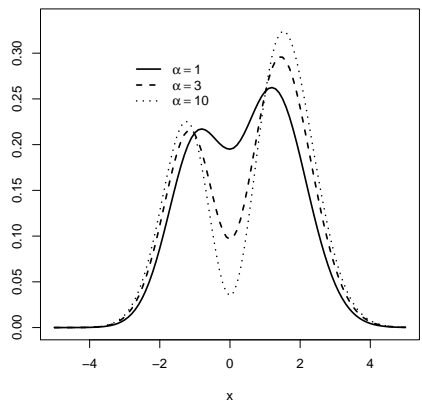,height=\figh in,width=\figw in}}& 
\subfigure[][$\gamma=0.9$]
   {\epsfig{file=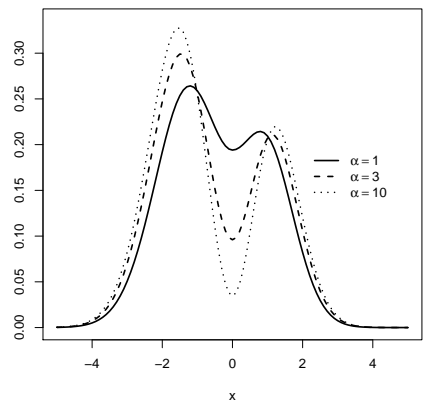,height=\figh in,width=\figw in}}   
\end{tabular}
\caption{Bimodal skew normal densities fixing the value of $\gamma$
  and varying $\alpha\in\{1,3,10\}$.}
\label{fig:bsn}
\end{figure}

\noindent For example, choosing $f(x)=\phi(x)$ in (\ref{eq:skew}), 
i.e. the density of a
standard normal distribution we obtain the bimodal skew normal
distribution with parameters $\alpha$ and $\gamma$ and
denote $X\sim BSN(\alpha,\gamma)$. This 
density is given by,
\begin{eqnarray}\label{eq:bsn}
s(x|\alpha,\gamma) 
&=& 
\left(\frac{1+\alpha x^2}{1+\alpha b_{\gamma}}\right)
\left(\frac{2}{\pi}\right)^{1/2}\frac{1}{(\gamma+1/\gamma)}\nonumber\\
&&
\exp\left\{-\frac{x^2}{2}
\left(
  \frac{1}{\gamma^2}I_{[0,\infty)}(x)+\gamma^2I_{(-\infty,0)}(x)
\right)\right\}, ~~x\in\mathbb{R},
\end{eqnarray}
and is depicted in Figure \ref{fig:bsn} for varying 
$\alpha\in\{1,3,10\}$ and fixing $\gamma>1$ (left panels) and
$\gamma<1$ (right panels).
For fixed $\alpha$, the position of the higher mode is controlled by
$\gamma$. As $\gamma>1$ (right skewness) increases density values are higher in
the right mode than in the left one as the original (unimodal) skewed
density puts more probability mass above zero. Actually, the left mode
is pushed towards zero as $\gamma$ increases above one. Of course the
reverse behaviour is observed when $0<\gamma<1$ (left skewness) decreases.


\begin{figure}[p]\centering
\begin{tabular}{cc}
\subfigure[][$\gamma=1.5$]
   {\epsfig{file=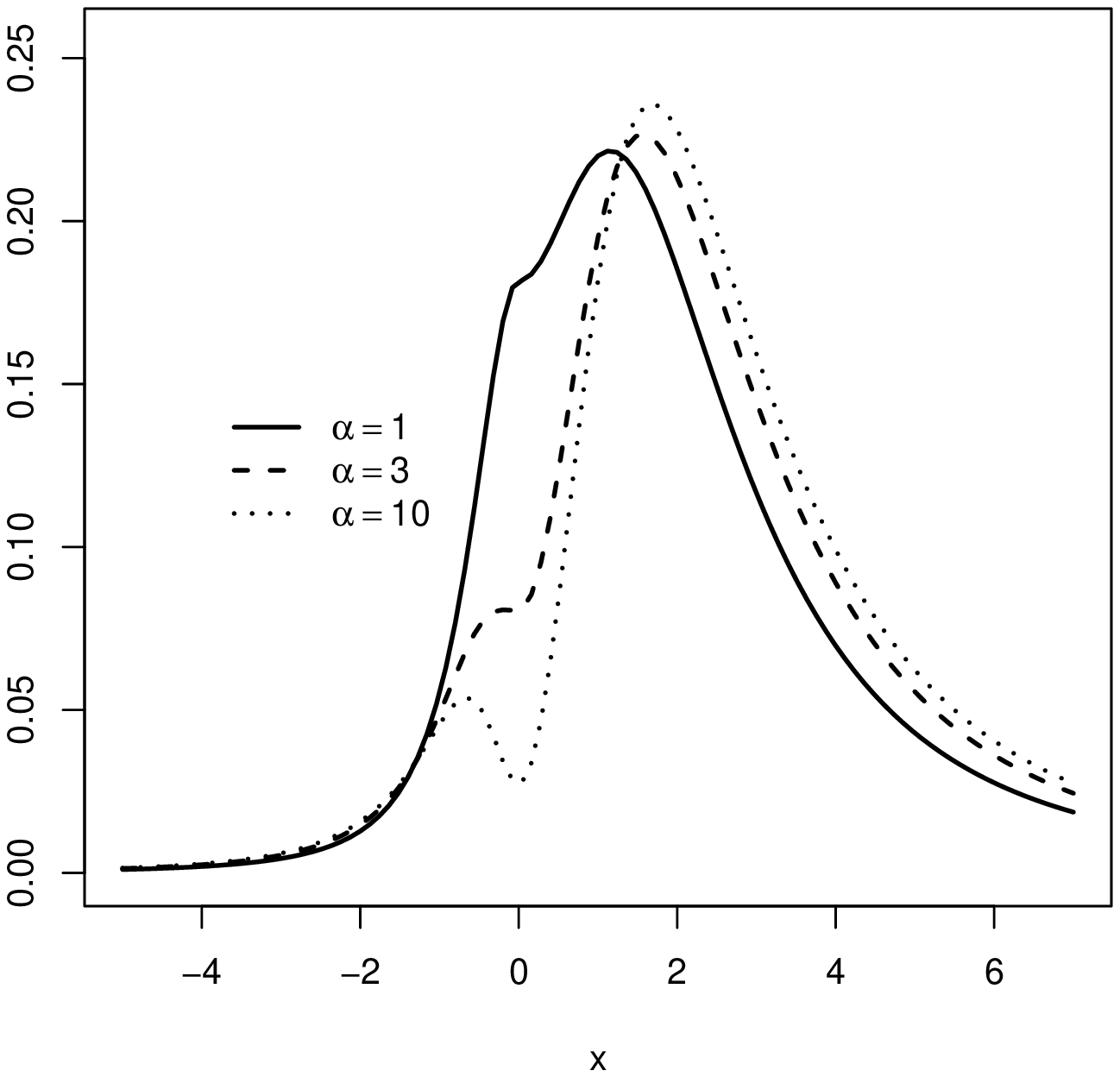,height=\figh in,width=\figw in}}& 
\subfigure[][$\gamma=0.5$]
   {\epsfig{file=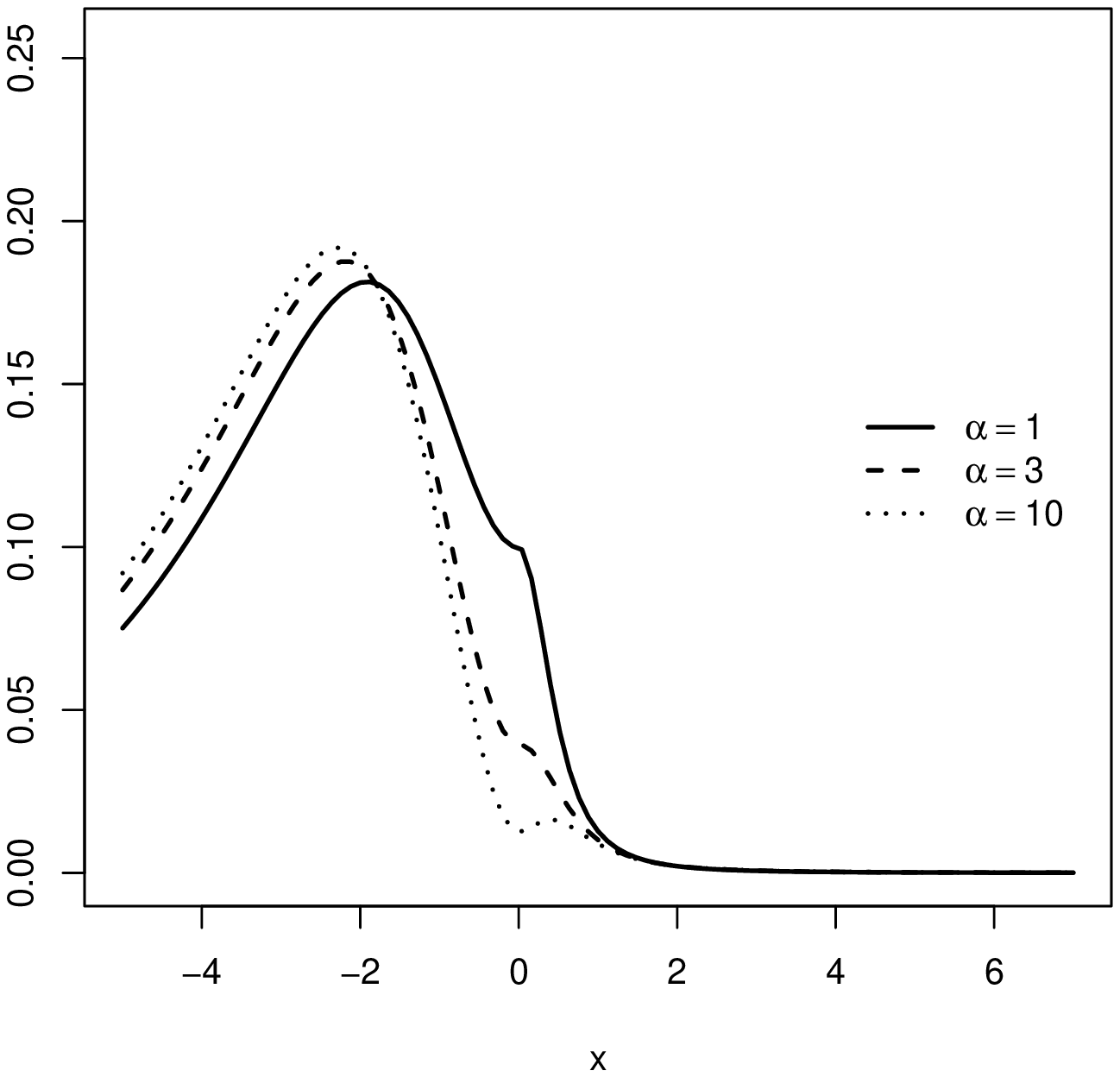,height=\figh in,width=\figw in}}\\
\subfigure[][$\gamma=1.1$]
   {\epsfig{file=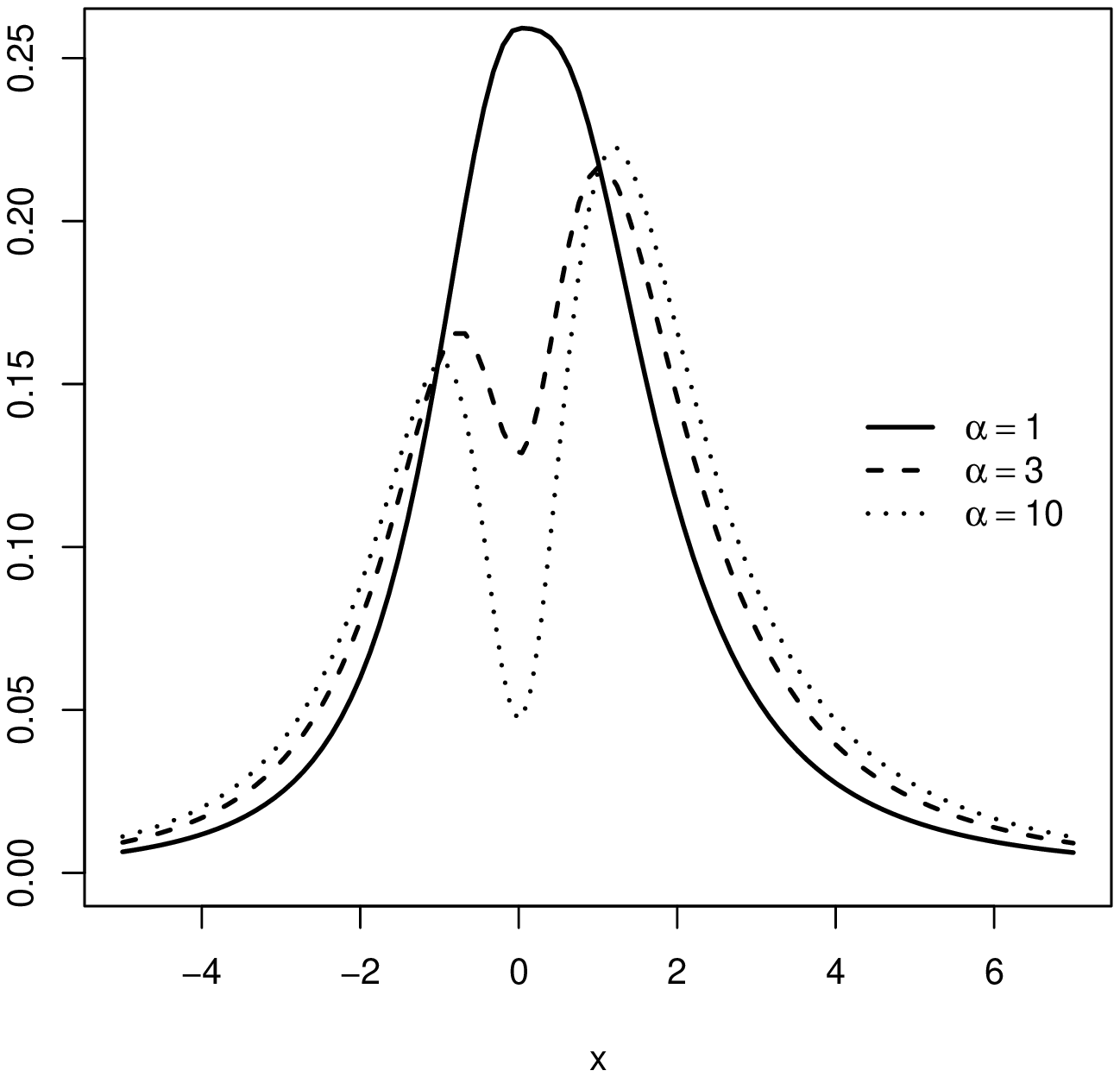,height=\figh in,width=\figw in}}& 
\subfigure[][$\gamma=0.9$]
   {\epsfig{file=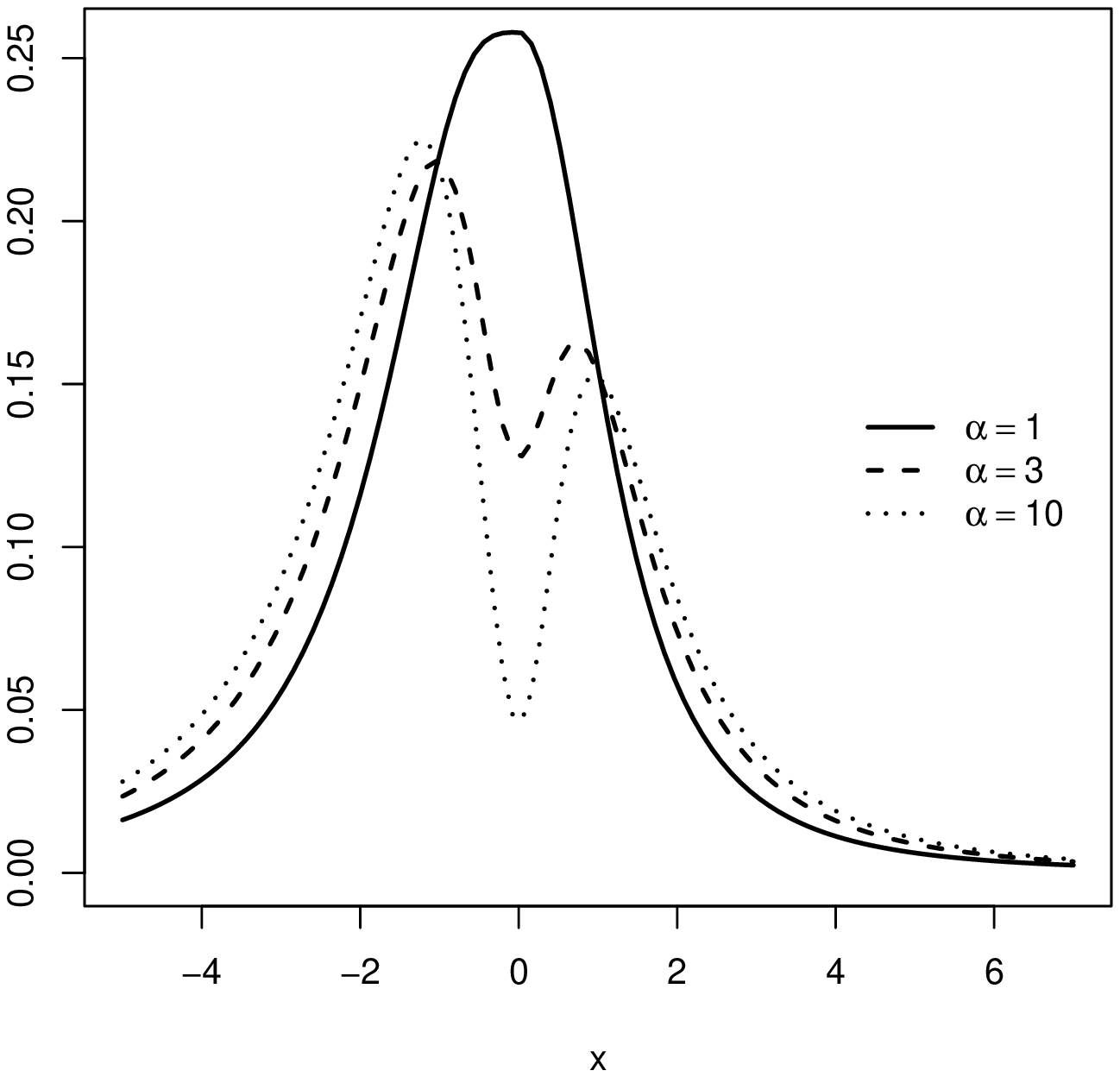,height=\figh in,width=\figw in}}   
\end{tabular}
\caption{Bimodal skew $t$ densities with $\nu=4$ fixing the value of $\gamma$
  and varying $\alpha\in\{1,3,10\}$.}
\label{fig:bsstd}
\end{figure}

Since they assign low probabilities to rare events, the family of
distributions presented above will fail to fit data with heavy tails
and we need to consider alternatives.
Choosing $f(\cdot)$ to be the standardized Student $t$ density (mean
zero and variance one) we obtain
the bimodal skewed Student distribution with parameters $\alpha$,
$\gamma$ and $\nu$ denoted $BSSTD(\alpha,\gamma,\nu)$ and density
function given by, 
\begin{eqnarray}\label{eq:bsstd}
s(x|\alpha,\gamma,\nu) 
&=& 
\left(\frac{1+\alpha x^2}{1+\alpha b_{\gamma}}\right)
\frac{2\Gamma(\frac{\nu+1}{2})}{\Gamma(\frac{\nu}{2})(\gamma+1/\gamma)[\pi(\nu-2)]^{1/2}}
\nonumber\\
&&
\left[
  1 + \frac{x^2}{\nu-2}
  \left\{\frac{1}{\gamma^2}I_{[0,\infty)}(x)+\gamma^2I_{[-\infty,0)}(x)\right\}
\right]^{-\frac{\nu+1}{2}},
\end{eqnarray}
for $x\in\mathbb{R}$ and $\nu>2$. 
This density is depicted in Figure \ref{fig:bsstd} for $\nu=4$, fixing
the value of $\gamma$ and varying $\alpha\in\{1,3,10\}$. It is clear
that, compared to the BSN case, events far apart in the tails will
receive higher probabilities under this family.
Note also that using this standardized version of the symmetric $t$
distribution allows us to keep the same expression for $b_{\gamma}$ in
both densities (\ref{eq:bsn}) and (\ref{eq:bsstd}) and propose a scale
mixture representation as follows.

{\proposition\label{p2} A random variable 
  $X\sim BSSTD(\alpha,\gamma,\nu)$ admits a scale mixture
  representation of $BSN(\alpha,\gamma)$ distributions with scale
  $\lambda^{-1/2}$ and mixing distribution $\lambda\sim Gamma(\nu/2,(\nu-2)/2)$.
}  
\vskip 0.5cm

\noindent {\bf Proof.}
Let $X|\lambda\sim BSN(\alpha,\gamma)$ with scale $\lambda^{-1/2}$ and density
given by,
\begin{eqnarray*}
s(x|\alpha,\gamma,\lambda) 
&=& 
\left(\frac{1+\alpha x^2}{1+\alpha b_{\gamma}}\right)
\left(\frac{2}{\pi}\right)^{1/2}\frac{\lambda^{1/2}}{(\gamma+1/\gamma)}\\
&&
\exp\left\{-\frac{\lambda x^2}{2}
\left(
  \frac{1}{\gamma^2}I_{[0,\infty)}(x)+\gamma^2I_{(-\infty,0)}(x)
\right)\right\}, ~~x\in\mathbb{R},
\end{eqnarray*}
and
$$
f(\lambda)= \frac{[(\nu-2)/2]^{\nu/2}}{\Gamma(\nu/2)}
\lambda^{\nu/2-1}\exp(-\lambda(\nu-2)/2).
$$
So, the marginal density of $X$ is given by,
\begin{eqnarray*}
s(x|\alpha,\gamma,\nu)
&=&
\left(\frac{1+\alpha x^2}{1+\alpha b_{\gamma}}\right)
\left(\frac{2}{\pi}\right)^{1/2}
\frac{1}{(\gamma+1/\gamma)}~
\frac{[(\nu-2)/2]^{\nu/2}}{\Gamma(\nu/2)}\times\\
&&
\int_0^{\infty}\lambda^{(\nu+1)/2-1}
\exp\left\{
-\frac{\lambda}{2}\left[(\nu-2)+x^2\gamma^{-2sign(x)}\right]
\right\} d\lambda\\
&=&
\left(\frac{1+\alpha x^2}{1+\alpha b_{\gamma}}\right)
\frac{2\Gamma(\frac{\nu+1}{2})}{\Gamma(\frac{\nu}{2})(\gamma+1/\gamma)[\pi(\nu-2)]^{1/2}}\\
&&
\left[1 + \frac{x^2}{\nu-2}\gamma^{-2sign(x)}\right]^{-\frac{\nu+1}{2}}
\end{eqnarray*}\qed

\noindent This representation will enable more efficient Bayesian
estimation via Markov chain Monte Carlo (MCMC)
algorithms using a data augmentation approach. As a by-product, the mixing
parameter $\lambda$ can be used to identify possible outliers.

{\proposition\label{p3} A random variable $X$ with a skewed normal
  distribution with scale $\lambda^{-1/2}$ admits the following
  hierarchical form,
\begin{eqnarray*}
  X|\gamma,\lambda,u &\sim&  
  SU\left(-\lambda^{-1/2}u^{1/2},\lambda^{-1/2}u^{1/2},\gamma\right)\\
  U &\sim& Gamma\left(\frac{3}{2}, \frac{1}{2}\right)
\end{eqnarray*}
where $SU(a,b,\gamma)$ denotes the skewed version of the
Uniform distribution on $(a,b)$. 
}
\vskip 0.5cm

\noindent {\bf Proof.} The density of the skewed version of a uniform
distribution on $(-\lambda^{-1/2}u^{1/2},\lambda^{-1/2}u^{1/2})$ is given by,

\begin{eqnarray*}
s(x|\gamma,\lambda,u) 
&=&
\frac{(\lambda/u)^{1/2}}{\gamma+1/\gamma}
\left[I(0<x<(u/\lambda)^{1/2}\gamma)+I(-(u/\lambda)^{1/2}/\gamma<x<0)\right]\\
&=&
\frac{(\lambda/u)^{1/2}}{\gamma+1/\gamma}
\left[I(u>\delta_1)I(x\ge 0) + I(u>\delta_2)I(x < 0)\right]
\end{eqnarray*}
where $\delta_1=\lambda x^2/\gamma^2$ and 
$\delta_2=\lambda x^2\gamma^2$. Now, integrating with respect to $u$ this
density times the density function of $u$ we obtain,

\begin{eqnarray*}
s(x|\gamma,\lambda) 
&=& 
\frac{\lambda^{1/2}}{\gamma+1/\gamma}~
\frac{(1/2)^{3/2}}{\Gamma(3/2)}\\
&&
\left[
  \int_{\delta_1}^{\infty}\exp(-u/2)du I_{[0,\infty)}(x)+
  \int_{\delta_2}^{\infty}\exp(-u/2)du I_{(-\infty,0)}(x)
\right]\\
&=&
\frac{\lambda^{1/2}}{\gamma+1/\gamma}~\left(\frac{2}{\pi}\right)^{1/2}
\exp\left\{
-\frac{\lambda x^2}{2}
\left[
  \frac{1}{\gamma^2}I_{[0,\infty)}(x)+\gamma^2I_{(-\infty,0)}(x)
\right]\right\}
\end{eqnarray*}
\qed
\vskip .5cm

\noindent Propositions \ref{p2} and \ref{p3} allow us to rewrite density
(\ref{eq:bsstd}) as the following scale mixture,

\begin{eqnarray*}
s(x|\alpha,\gamma,\nu) &=&
\int_0^{\infty}\int_0^{\infty}
\left(\frac{1+\alpha x^2}{1+\alpha b_{\gamma}}\right)
s(x|\gamma,\lambda,u)\times\\
&&
f_{{\tiny G}}\left(u\bigg|\frac{3}{2},\frac{1}{2}\right)
f_{{\tiny G}}\left(\lambda\bigg|\frac{\nu}{2},\frac{\nu-2}{2}\right)dud\lambda.
\end{eqnarray*}
where
$f_{{\tiny G}}(\cdot|a,b)$ denotes the density of a Gamma distributed random
variable with mean $a/b$ and variance $a/b^2$.

\section{A Wider Class of Distributions}

\begin{figure}[p]\centering
\begin{tabular}{cc}
\subfigure[][$\gamma=1.5$]
   {\epsfig{file=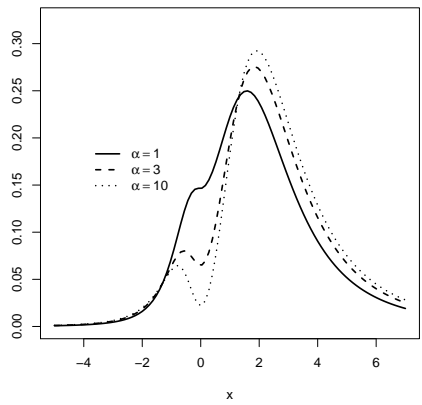,height=\figh in,width=\figw in}}& 
\subfigure[][$\gamma=0.5$]
   {\epsfig{file=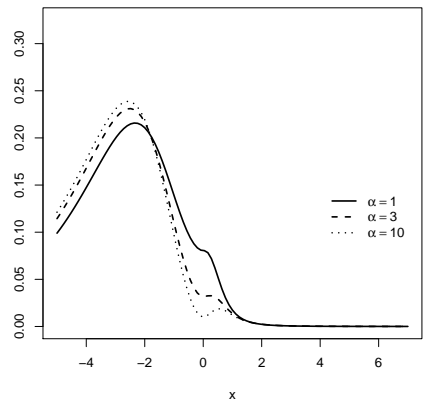,height=\figh in,width=\figw in}}\\
\subfigure[][$\gamma=1.1$]
   {\epsfig{file=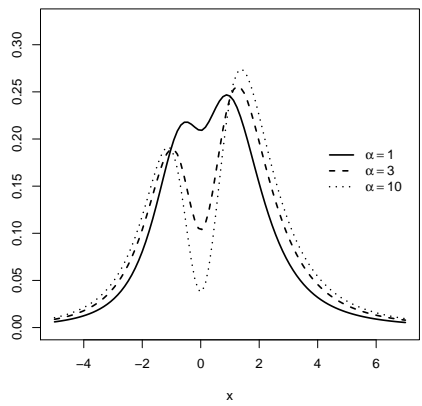,height=\figh in,width=\figw in}}& 
\subfigure[][$\gamma=0.9$]
   {\epsfig{file=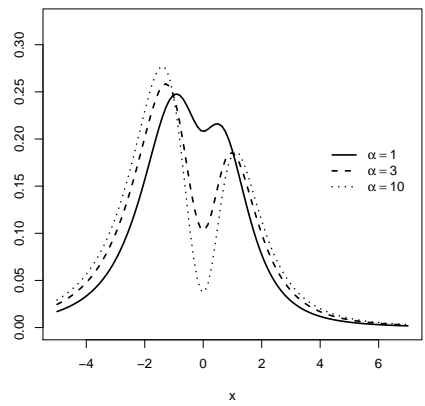,height=\figh in,width=\figw in}}   
\end{tabular}
\caption{Bimodal skew generalized $t$ densities with $p=2.3$, $q=2$,
  fixing the value of $\gamma$ and varying $\alpha\in\{1,3,10\}$.}
\label{fig:bsgt}
\end{figure}

\begin{figure}[p]\centering
\begin{tabular}{cc}
\subfigure[][$\gamma=1.5$]
   {\epsfig{file=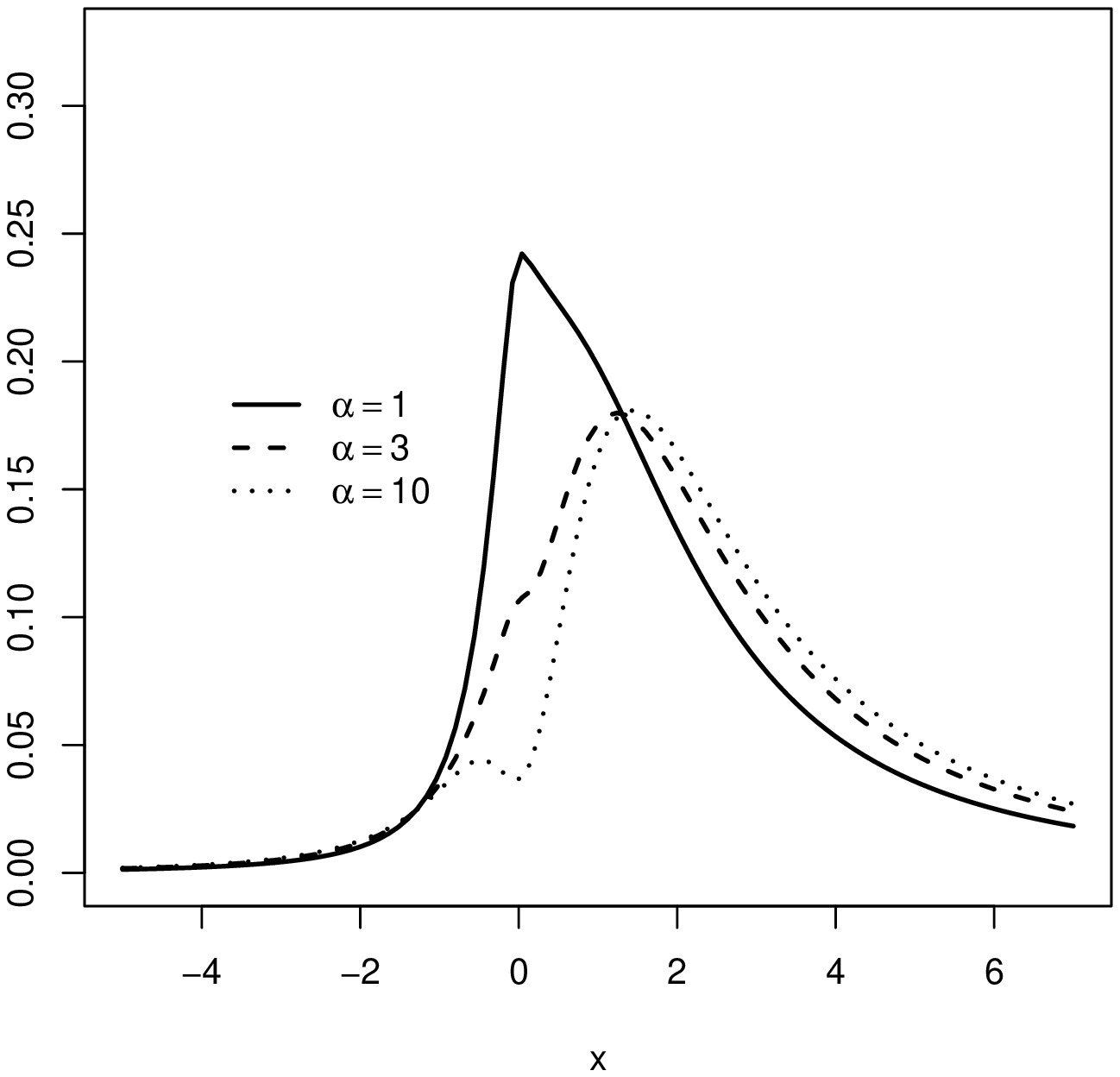,height=\figh in,width=\figw in}}& 
\subfigure[][$\gamma=0.5$]
   {\epsfig{file=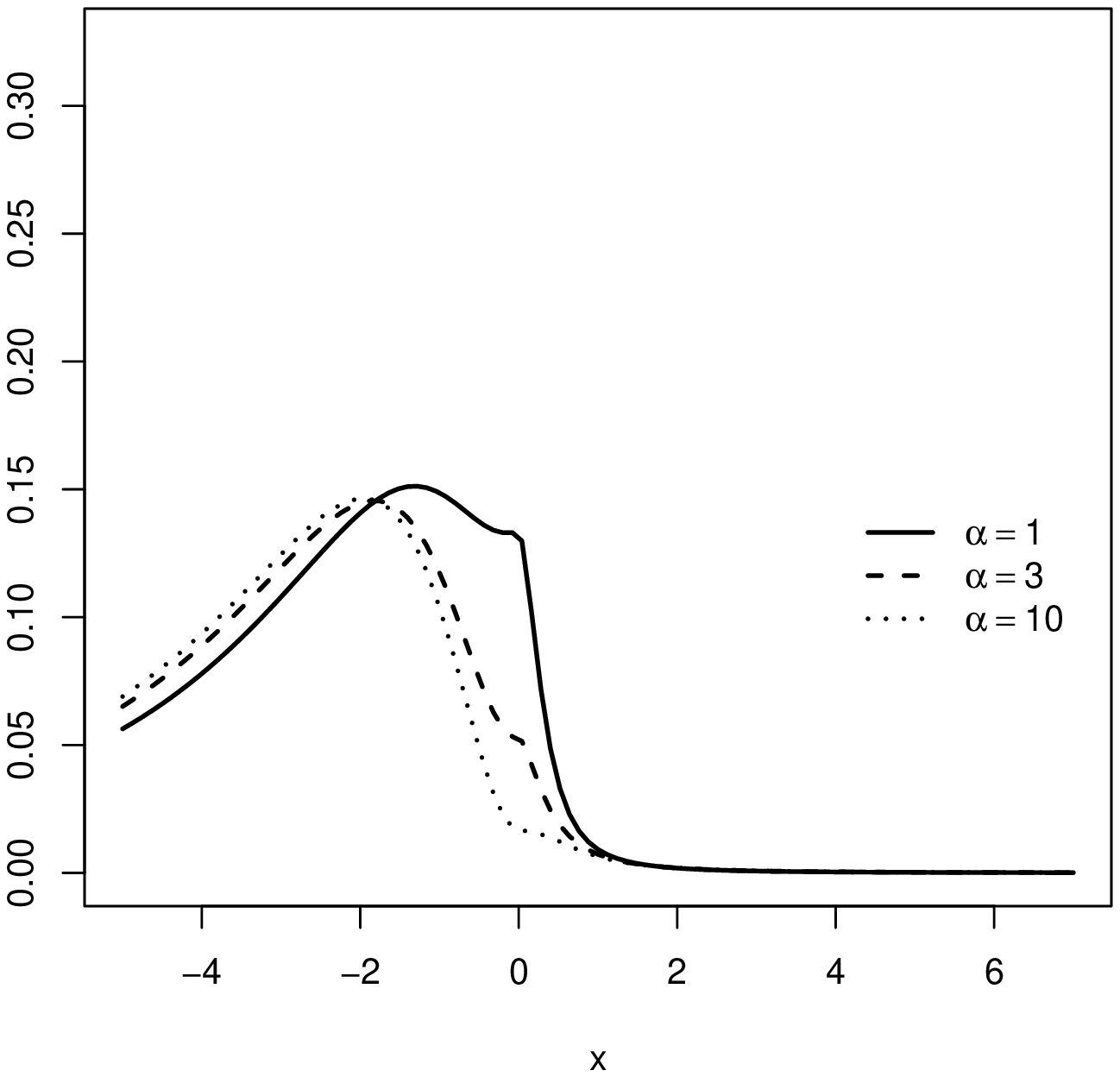,height=\figh in,width=\figw in}}\\
\subfigure[][$\gamma=1.1$]
   {\epsfig{file=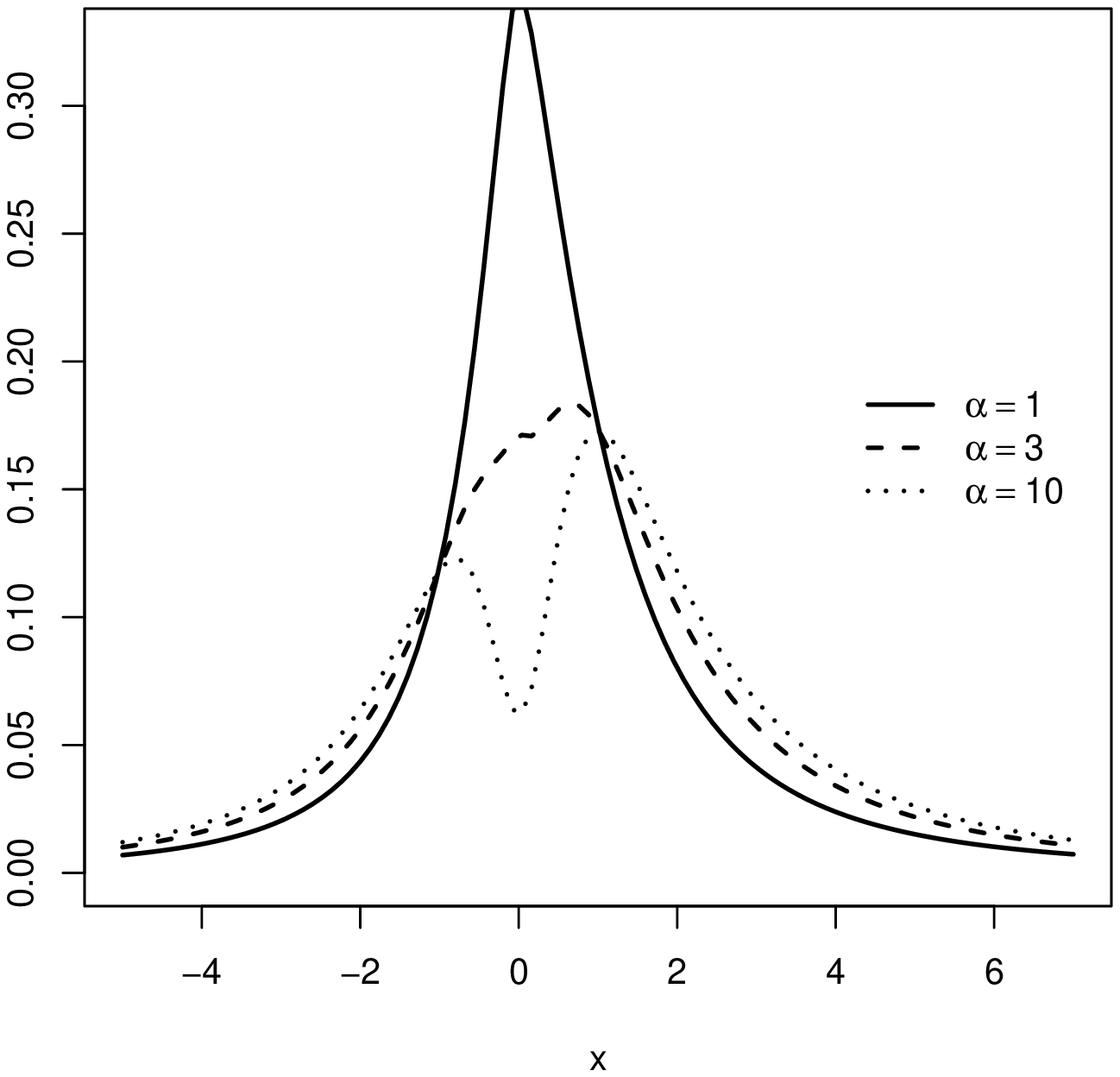,height=\figh in,width=\figw in}}& 
\subfigure[][$\gamma=0.9$]
   {\epsfig{file=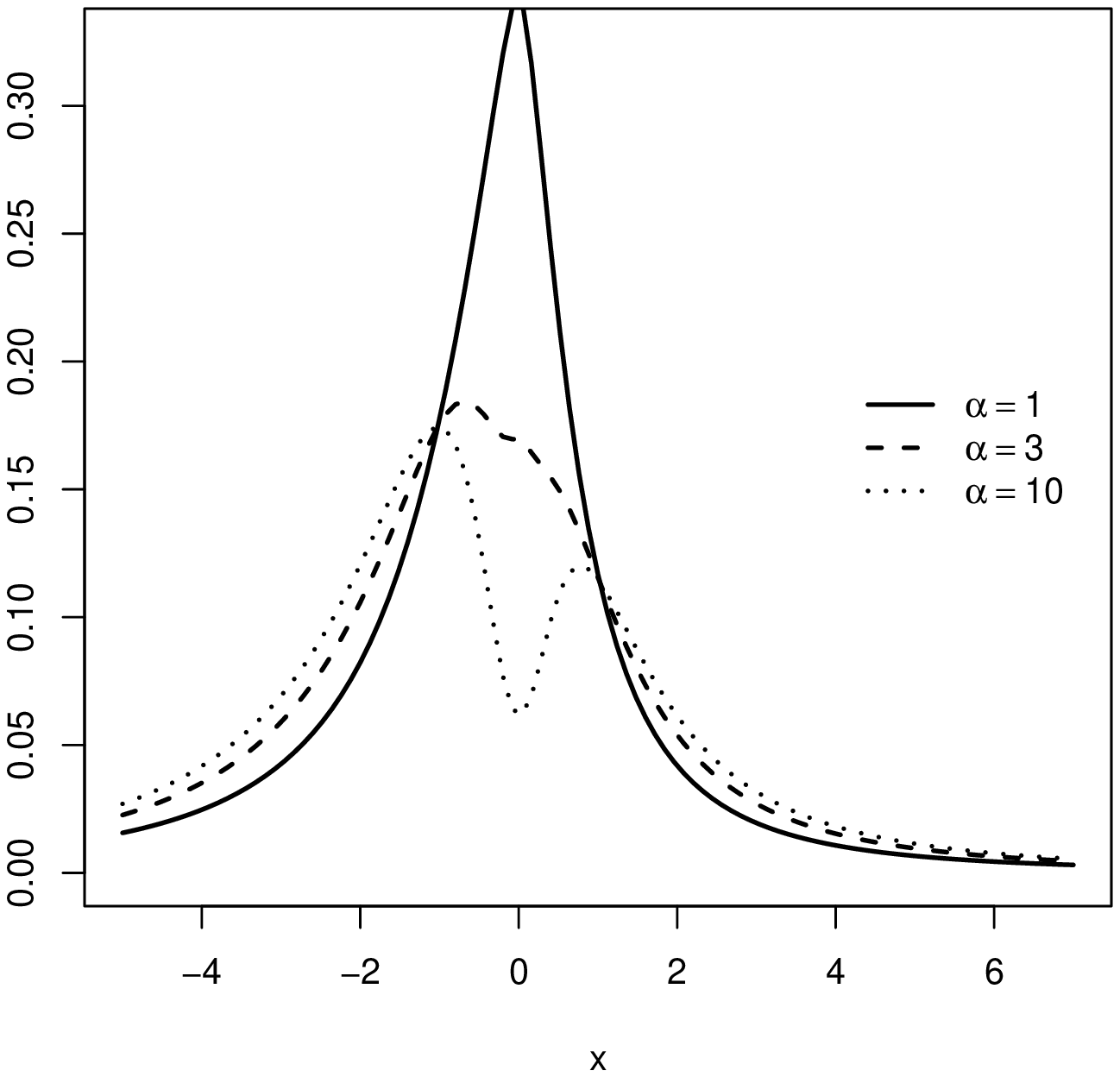,height=\figh in,width=\figw in}}   
\end{tabular}
\caption{Bimodal skew generalized $t$ densities with $p=1.7$,
  $q=2$, fixing the value of $\gamma$ and varying $\alpha\in\{1,3,10\}$.}
\label{fig:bsgt1}
\end{figure}

\citeN{mcdonald-newey88} introduced a flexible symmetric and unimodal
distribution as another robust alternative to the normal distribution
which they called the generalized $t$ distribution.
Its density function with location zero and scale one is given by,
\begin{eqnarray}\label{eq:gt}
f(x) &=&
\frac{p\Gamma\left(q+\frac{1}{p}\right)}
{2q^{1/p}\Gamma\left(\frac{1}{p}\right)\Gamma\left(q\right)}
\left(1+\frac{1}{q}\left|x\right|^p\right)^{-(q+1/p)}\nonumber\\
&=&
\frac{p}{2q^{1/p}B(1/p,q)}
\left(1+\frac{1}{q}\left|x\right|^p\right)^{-(q+1/p)}
~x\in\mathbb{R}
\end{eqnarray}
where $B(a,b)=\Gamma(a)\Gamma(b)/\Gamma(a+b)$ is the Beta function,
$p>0$ and $q>0$ are two shape parameters. We refer to
this distribution as $GT(p,q)$. Larger values of $p$ and $q$ yield a
density with thinner tails than the 
normal while smaller values are associated with thicker tailed
densities. Also, it includes other well known symmetric unimodal
distributions as special or limiting cases. In particular, the
variance exists when $pq>2$ and is given by,
$$
Var(X)=\frac{q^{2/p}\Gamma(3/p)\Gamma(q-2/p)}{\Gamma(1/p)\Gamma(q)}=
\frac{q^{2/p}B(3/2,q-2/p)}{B(1/p,q)}.
$$
Therefore, the standardized version of density (\ref{eq:gt}) is given by,
\begin{equation}
f(x) = 
\frac{p}{2\delta q^{1/p}B(1/p,q)}
\left(1+\frac{1}{q}\left|\frac{x}{\delta}\right|^p\right)^{-(q+1/p)}
~x\in\mathbb{R}
\end{equation}
where
$$
\delta=
\left(\frac{q^{2/p}B(3/2,q-2/p)}{B(1/p,q)}\right)^{-1/2}.
$$
Using this standardized version of the generalized $t$ distribution we obtain
the bimodal skewed generalized $t$ distribution with parameters
$\alpha$, $\gamma$, $p$ and $q$ denoted $BSGT(\alpha,\gamma,p,q)$ and
density given by,
\begin{eqnarray}\label{eq:bsgt}
s(x|\alpha,\gamma,p,q) 
&=&
\left(\frac{1+\alpha x^2}{1+\alpha b_{\gamma}}\right)
\frac{p}{\delta (\gamma+1/\gamma) q^{1/p}B(1/p,q)}\nonumber\\
&&
\left[1+\frac{1}{q}\left|\frac{x}{\delta}\right|^p
\left\{
\frac{1}{\gamma^{p}}I_{[0,\infty)}(x)+\gamma^{p}I_{(-\infty,0)}(x)
\right\}\right]^{-(q+1/p)}.
\end{eqnarray}
Again using a standardized version of the original symmetric
distribution allows us to keep the same expression for the second
moment $b_{\gamma}$. It is not difficult to see that setting $p=2$ we
recover the bimodal skewed $t$ distribution with tail parameter
$\nu=2q$. The bimodal skewed normal is then obtained when $p=2$ and
$q\rightarrow\infty$. Density (\ref{eq:bsgt}) is depicted in Figure
\ref{fig:bsgt} with $p=2.3$, $q=2$, fixing the value of $\gamma$
and varying $\alpha\in\{1,3,10\}$. Parameter $p$ has a larger
influence on the shape of the density than $q$, a feature inhereted
from the symmetric version of the GT distribution. This is
illustrated in Figure \ref{fig:bsgt1} where we set $p=1.7$

One feature of the symmetric generalized $t$ distribution is that it
can be represented as a scale mixture of an exponential power
distribution (\citeNP{boxtiao73}) with a generalized Gamma as the mixing
distribution, a result obtained by \citeN{arslang03}. In the next
proposition we extend this representation to the bimodal skewed GT
distribution with density (\ref{eq:bsgt}).

{\proposition\label{p4} A random variable $X\sim
  BSGT(\alpha,\gamma,p,q)$ with density (\ref{eq:bsgt})
  admits a scale mixture representation with the following
  hierarchical form,
  \begin{eqnarray*}
    X|S=s &\sim& BSEP\left(2^{-1/p}s^{-1/2}\left[\frac{\Gamma(q-2/p)}{\Gamma(q)}\right]^{-1/2}, ~p\right)\\
      S   &\sim& GG\left(\frac{p}{2},1,q\right),
  \end{eqnarray*}
  where $BSEP(\lambda,p)$ denotes a bimodal skewed exponential power
  distribution with scale $\lambda$ and tail parameter $p$ and
  $GG(\cdot,\cdot,\cdot)$ denotes the generalized Gamma distribution}.
\vskip .5cm

\noindent {\bf Proof.} The density of a (symmetric) standardized exponential power
distribution with tail parameter $p$ is given by,
\begin{eqnarray*}
f(x|p) &=& 
\frac{p}{\phi\Gamma(1/p) ~2^{1+1/p}}
\exp\left\{-\frac{1}{2}\left|\frac{x}{\phi}\right|^p
\right\}
\end{eqnarray*}
where $\phi = [2^{2/p}\Gamma(3/p)/\Gamma(1/p)]^{-1/2}$.
The skewed version of this exponential power
distribution with scale given in the proposition,
tail parameter $p$ and noting that 
$$
\phi=2^{1/p} q^{-1/p} \left[\frac{\Gamma(q-2/p)}{\Gamma(q)}\right]^{1/2}\delta
$$
is then given by,
\begin{eqnarray*}
s(x|p,q,\gamma,s) &=& 
\frac{2^{1/p}s^{1/2} ~p}{\delta ~q^{1/p} ~\Gamma(1/p) ~2^{1+1/p}}~
\frac{2}{\gamma+1/\gamma}\\
&&
\exp\left\{-\frac{1}{2}\left|\frac{x}{2^{-1/p}s^{-1/2}q^{1/p}\delta}\right|^p
\left(\frac{1}{\gamma^p}I_{[0,\infty)}(x)+\gamma^pI_{(-\infty,0)}(x)\right)
\right\}
\end{eqnarray*}
while the density of a generalized Gamma distribution with parameters
$1/p$, 1 and $q$ is given by,
$$
f(s) = \frac{p}{2\Gamma(q)}s^{pq/2-1}\exp(-s^{p/2}).
$$
Since the original symmetric density is in its standardized form we
have the same expression for the second moment $b_{\gamma}$ of the
skewed density. It then follows that,

\begin{eqnarray*}
s(x|\alpha,\gamma,p,q) 
&=& \int_0^{\infty}
\left(\frac{1+\alpha x^2}{1+\alpha b_{\gamma}}\right)
s(x|p,q,\gamma,s)f(s) ds\\
&=&
\left(\frac{1+\alpha x^2}{1+\alpha b_{\gamma}}\right)
\frac{p^2}{2~\delta~q^{1/p}\Gamma(1/p)\Gamma(q)}~
\frac{1}{\gamma+1/\gamma}\times\\
&&
\int_0^{\infty} s^{pq/2+1/2-1}
\exp\left\{-s^{p/2}\left[1+\frac{1}{q}
  \left|\frac{x}{\delta}\right|^p (\gamma^{-p})^{sign(x)}
  \right]\right\} ds.
\end{eqnarray*}
Now, setting $y=s^{p/2}$ the last integral is rewritten as,
$$
\frac{2}{p}\int_0^{\infty} y^{(q+1/p)-1}
\exp\left\{-y\left[1+\frac{1}{q}
\left|\frac{x}{\delta}\right|^p (\gamma^{-p})^{sign(x)}
\right]\right\} dy,
$$
and finally,
\begin{eqnarray*}
s(x|\alpha,\gamma,p,q) &=&
\left(\frac{1+\alpha x^2}{1+\alpha b_{\gamma}}\right)
\frac{p\Gamma(q+1/p)}{\delta q^{1/p}(\gamma+1/\gamma)\Gamma(1/p)\Gamma(q)}\\
&&
\left[1+\frac{1}{q}
\left|\frac{x}{\delta}\right|^p (\gamma^{-p})^{sign(x)}
\right]^{-(q+1/p)}.
\end{eqnarray*}
\qed

In what follows we propose an alternative representation
based on the skewed version of the uniform distribution used in
Proposition \ref{p3}. \citeN{choy-chan} had already proposed an
alternative representation for the symmetric generalized $t$
distribution based on
a scale mixture of (symmetric) uniform distributions. This was latter
extended in \citeN{ehl2015d} for a skewed version.

{\proposition \label{p5} 
A random variable $X\sim BSGT(\alpha,\gamma,p,q)$ with density
(\ref{eq:bsgt}) admits a  
scale mixture representation with the following hierarchical form,
\begin{eqnarray*}
  X|u,s,\alpha,\gamma &\sim& BSU
  \left(-a,a,\alpha,\gamma\right)\\
  U &\sim& Gamma\left(1+\frac{1}{p},1\right)\\
  S &\sim& GG\left(\frac{p}{2},1,q\right)
\end{eqnarray*}
where
$$
a=2^{-1/p}s^{-1/2}\left[\frac{\Gamma(q-2/p)}{\Gamma(q)}\right]^{-1/2}u^{1/p}
$$
}

\section{Bayesian Inference}

Following \citeN{fsteel98}, we shall use a Gamma($a,b$) prior
distribution on
$\phi=\gamma^2$ which is the ratio of probability masses above and below
the mode, i.e. $\phi=\gamma^2=Pr(X\ge 0)/Pr(X < 0)$. 
For observed data $\bfx=(x_1,\dots,x_n)$ the likelihood function 
in the bimodal skewed normal model is given by
\begin{eqnarray*}
s(\bfx|\alpha,\phi) &\propto&
(1+\alpha b_{\phi})^{-n}
\phi^{n/2}(1+\phi)^{-n}
\prod_{i=1}^n (1+\alpha x_i^2)\\
&&
\exp\left\{-\frac{1}{2}\sum_{i=1}^n x_i^2
\left(
  \frac{1}{\gamma^2}I_{[0,\infty)}(x_i)+\gamma^2I_{(-\infty,0)}(x_i)
\right)\right\}
\end{eqnarray*}
where 
$$
b_{\phi} = \frac{1+\phi^3}{\phi(1+\phi)}.
$$

\noindent The complete conditional distributions of $\phi$ and $\alpha$
are then given by,
\begin{eqnarray*}
f(\phi|\bfx,\alpha) 
&\propto&
(1+\alpha b_{\phi})^{-n}\phi^{a+n/2-1}(1+\phi)^{-n}\\
&&
\exp\left\{-\frac{1}{2}\sum_{i=1}^n x_i^2\phi^{-sign(x_i)}-b\phi\right\}.
\end{eqnarray*}

\begin{eqnarray*}
f(\alpha|\bfx,\phi) 
&\propto&
(1+\alpha b_{\phi})^{-n}\prod_{i=1}^n (1+\alpha x_i^2) ~f(\alpha).
\end{eqnarray*}

If we now assume a bimodal skew $t$ distribution we need to assign a prior
distribution for the tail parameter $\nu$. Here, we follow
\citeN{deschamps06} and use a translated exponential distribution with density,
$$
p(\nu)=\beta \exp\{-\beta(\nu-2)\} I(\nu>2).
$$

\noindent Using the scale mixture representation, each observation
$X_i$ is associated with the mixing parameter $\lambda_i$ and we
assume that they are a priori independent. The
complete conditional densities are given by, 

\begin{eqnarray*}
f(\alpha|\bfx,\bflambda,\gamma,\nu)
&\propto&
s(\bfx|\alpha,\gamma,\bflambda) ~f(\alpha)\\
&\propto&
(1+\alpha b_{\phi})^{-n}\prod_{i=1}^n (1+\alpha x_i^2) ~f(\alpha).
\end{eqnarray*}

\begin{eqnarray*}
f(\phi|\bfx,\bflambda,\alpha,\nu)
&\propto&
s(\bfx|\alpha,\gamma,\bflambda) ~f(\phi)\\\\
&\propto&
(1+\alpha b_{\phi})^{-n}\phi^{a+n/2-1}(1+\phi)^{-n}\\
&&
\exp\left\{
-\frac{1}{2}\sum_{i=1}^n \lambda_ix_i^2\phi^{-sign(x_i)}-b\phi
\right\}.
\end{eqnarray*}

\begin{eqnarray*}
f(\nu|\bfx,\bflambda,\alpha,\gamma)
&\propto&
f(\nu) \prod_{i=1}^n f(\lambda_i|\nu)\\
&\propto&
\exp\{-\beta(\nu-2)\}\\
&&
\prod_{i=1}^n \frac{[(\nu-2)/2]^{\nu/2}}{\Gamma(\nu/2)}
\lambda_i^{(\nu-2)/2}\exp\left(-\frac{\lambda_i(\nu-2)}{2}\right)\\
&\propto&
\frac{[(\nu-2)/2]^{n\nu/2}}{\Gamma^n(\nu/2)}
\exp\left\{-\nu\left(
\beta+\frac{1}{2}\sum_{i=1}^n(\lambda_i-\log\lambda_i)
\right)\right\}.
\end{eqnarray*}

\begin{eqnarray*}
  f(\bflambda|\bfx,\alpha,\gamma,\nu)&\propto&
  \prod_{i=1}^n s(x_i|\alpha,\gamma,\lambda_i) f(\lambda_i|\nu)\\
  &\propto&
  \prod_{i=1}^n \lambda_i^{(\nu+1)/2-1}
  \exp\left\{-\frac{\lambda_i}{2}(\nu-2+x_i^2\gamma^{-2sign(x_i)})\right\}
\end{eqnarray*}
so the complete conditional distribution of each mixing parameter is
given by,
$$
\lambda_i|\bfx,\bflambda_{-i},\alpha,\gamma,\nu\sim 
Gamma\left(\frac{\nu+1}{2},\frac{\nu-2+x_i^2\gamma^{-2sign(x_i)}}{2}\right),
$$
being easily sampled from. However, the complete conditional
distributions of $\alpha$, $\phi$ and $\nu$ are not of any standard form.

%

\end{document}

%% file: mymaths.tex
\newcommand{\bflambda }{\hbox{\boldmath$\lambda$}}

\newcommand{\bfx}{\hbox{\boldmath$x$}}


%% file: theorems.tex
\newtheorem{proposition}{Proposition}[section]

\newcounter{qqq}[section]
\vspace{0.5cm}
\vspace{0.5cm}